\documentclass[12pt]{article}

\usepackage[intlimits,tbtags]{amsmath}
\usepackage{amsfonts,amsthm}
\usepackage{latexsym}

\newcommand{\E}{\mathrm{e}}
\newcommand{\I}{\mathrm{i}}

\newcommand{\tr}{\triangleright}
\newcommand{\act}{\triangleright}
\newcommand{\R}{\mathcal{R}}
\newcommand{\RI}{ \mathcal{R}_{\mathrm{I}}{} }
\newcommand{\RII}{ \mathcal{R}_{\mathrm{II}}{} }
\newcommand{\id}{\mathrm{id}}
\newcommand{\op}{\mathrm{op}}
\newcommand{\Xcal}{\mathcal{X}}
\newcommand{\h}{\hbar}

\newcommand{\Xsquare}{{x^2}}
\newcommand{\qNum}[1]{{(#1)_{q^2}}}

\newcommand{\Uqsu}{{\mathcal{U}_q(\mathrm{su}_2)}}
\newcommand{\Usu}{{\mathcal{U}(\mathrm{su}_2)}}

\newcommand{\suq}{{\mathcal{U}_q(\mathrm{su}_2)}}
\newcommand{\slC}{{\mathcal{U}_q(\mathrm{sl}_2(\mathbb{C})) }}

\newcommand{\SUq}{{SU_q(2)}}

\newcommand{\Mink}{\Xcal}

\newtheorem{Definition}{Definition} 

\newtheorem{Proposition}{Proposition}

\newcommand{\Lind}{{3\!\setminus\! 0}}
\newcommand{\Mat}[1]{{#1}}

\begin{document}

\vspace{3em}
\begin{center}
  
  {\Large{\bf Separation of noncommutative differential calculus on
      quantum Minkowski space}}

\vspace{3em}

\textbf{Fabian Bachmaier$^{1}$, Christian Blohmann$^{2,3}$}

\vspace{1em}
 
${}^1$Arnold Sommerfeld Center for Theoretical Physics\\
Universit\"at M\"unchen, Fakult\"at f\"ur Physik\\
Theresienstr.\ 37, 80333 M\"unchen, Germany\\[1em]


${}^2$Department of Mathematics, University of California\\
Berkeley, California 94720-3840\\[1em]

${}^3$International University Bremen, School of Engineering and Science, Campus Ring 1, 28759 Bremen, Germany
\\[1em]

\end{center}

\vspace{1em}

\begin{abstract}
Noncommutative differential calculus on quantum Minkowski space is not separated with respect to the standard generators, in the sense that partial derivatives of functions of a single generator can depend on all other generators. It is shown that this problem can be overcome by a separation of variables. We study the action of the universal $L$-matrix, appearing in the coproduct of partial derivatives, on generators. Powers of he resulting quantum Minkowski algebra valued matrices are calculated. This leads to a nonlinear coordinate transformation which essentially separates the calculus. A compact formula for general derivatives is obtained in form of a chain rule with partial Jackson derivatives. It is applied to the massive quantum Klein-Gordon equation by reducing it to an ordinary $q$-difference equation. The rest state solution can be expressed in terms of a product of $q$-exponential functions in the separated variables.
\end{abstract}

\newpage

\section{Introduction}

\subsection{Quantum theory on noncommutative space time}

Much of the development of quantum theories on noncommutative space-times was and still is driven by the question whether noncommutative geometry might lead to an ultra-violet regularization of quantum field theory, as it was suggested by Heisenberg as early as 1938 \cite{Heisenberg:1938}. For the simplest conceivable examples of noncommutative geometries, where the commutator of the coordinates is constant, these hopes for regularization were met with disappointment. Such theories exhibit a mixing of ultra-violet and infra-red cutoff scales \cite{Minwalla:1999,Matusis:2000}, which has been understood recently on the level of renormalization \cite{Grosse:2005}. For that reason these simple noncommutative spaces, while interesting objects on their own, hardly seem to improve the divergent behavior of quantum field theory at all. The natural next step is to turn again to more complicated noncommutative spaces, such as quantum spaces with Lie type or homogeneous commutation relations.

Quantum Minkowski space is one of the most realistic examples of a quantum deformation of space-time \cite{Carow-Watamura:1990} with a rich and fairly well understood mathematical structure: It is four dimensional in the sense that it is generated as algebra by four coordinates such that the ordered monomials are a Poincar\'e-Birkhoff-Witt basis. The time coordinate is central, which is important for a causal interpretation of quantum theory \cite{Doplicher:1995}.  By construction, it is a module algebra with respect to the quantum Poincar\'e algebra \cite{Ogievetskii:1992a}, so it has a well defined quantum symmetry structure \cite{Faddeev:1990}. The action of the generators of the inhomogeneous part of the quantum Poincar\'e algebra, the momenta, on quantum Minkowski space induces a first order covariant differential calculus \cite{Majid:1993}.
 
With these key mathematical structures present, much of the construction of quantum theory on commutative space-time was mimicked within an algebraic and representation theoretic approach. Sections of noncommutative space-time were given by spectra of noncommutative coordinates, which yield discrete space-time lattices \cite{Cerchiai:1998}. Free elementary particles were constructed as irreducible representations of the quantum Poincar\'e algebra, including Wigner representations of spin \cite{Blohmann:2001a}. The free particles were shown to obey quantum wave equations which project reducible quantum Lorentz spinor representations on their irreducible subrepresentations \cite{Blohmann:2001b}. The coupling of free fields to gauge fields was studied systematically within noncommutative gauge theory \cite{Madore:2000b,Jurco:2001}. Quantum statistics of the tensor representation of a second quantized free theory was studied within the framework of braided tensor categories \cite{Majid} using Drinfeld twists \cite{Fiore:1996}. Only a small selection of the numerous contributions to this program can be cited here. At the time of this writing, most aspects of free quantum theory on quantum Minkowski space are understood, and the construction of interaction terms on the level of single particle field equations is known. What is not known is how to solve these field equations by noncommutative wave functions, not even for the free case. However, this will be indispensable if free particles are to be coupled by multiplying their noncommutative wave functions within the algebra of quantum space-time.

We remark, that there are other approaches which involve wave equations on noncommutative space-time which are inequivalent to the ones considered here: In \cite{Dobrev:1994} covariance arguments have been used to find wave equations for a certain deformation of the conformal symmetry algebra, which does not contain the $q$-Poincar\'e algebra we consider here. In order to construct $q$-deformed relativistic wave equations various other methods have been proposed, based on $q$-Clifford algebras \cite{Schirrmacher:1992}, $q$-deformed co-spinors \cite{Pillin:1994b}, or various abstract differential calculi on quantum spaces \cite{Song:1992,Meyer:1995,Podles:1996}, each leading to different results.

\subsection{Noncommutative differential calculus}

The free wave equations on quantum Minkowski space can be determined uniquely on a purely representation theoretic level. To give the simplest example: Just as in the undeformed case, the momentum square $p_\mu p^\mu$ is central within the quantum Poincar\'e algebra. Therefore, it must be represented within the irreducible representation corresponding to a free particle by a multiple of the identity operator,
\begin{equation}
\label{eq:KleinGordon1}  
  p_\mu p^\mu = m^2 \,,
\end{equation}
$m^2$ being the square of the particle's mass. This line of reasoning is purely representation theoretic and does not depend in any way on how the momenta act on noncommutative wave functions, that is, on elements of quantum Minkowski space.

In order to interpret Eq.~\eqref{eq:KleinGordon1} as Klein-Gordon equation on noncommutative space-time we have to let the momenta act on noncommutative wave functions. Momenta then correspond to partial derivatives, $p^\mu = \I \partial^\mu$, defined by the noncommutative first order differential calculus. This yields a wave equation given by a noncommutative partial differential equation,
\begin{equation}
\label{eq:KleinGordon2}  
  \partial_\mu \partial^\mu \tr \psi = - m^2 \psi  \,,
\end{equation}
where $\psi$ is an element of the algebra of quantum Minkowski space. There is a considerable amount on literature on first order differential calculi on noncommutative algebras, mostly on quantum groups \cite{Woronowicz:1989} and to some extent on homogeneous quantum spaces \cite{Wess:1991,Schmuedgen}. Most mathematical work has concentrated on the construction, structural analysis, and classification of differential calculi. However, hardly any literature on the subject has dealt with concrete calculations within these differential calculi. On quantum Minkowski space a differential calculus has been known for some time, and can be deduced most elegantly from the coproduct of momenta \cite{Majid:1993}. The coproduct defines the action of partial derivatives on noncommutative functions recursively. While the recursion relations can be used to expand the derivatives in terms of nested summations and partition functions, as it was carried out in a very detailed manner in \cite{Bauer:2002}, this does not solve the main computational problem which arises with wave equations like Eq.~\eqref{eq:KleinGordon2}.

Trying to solve Eq.~\eqref{eq:KleinGordon2} by a brute force calculation with a general ansatz for the wave function ought to work in principle, but turns out to run into considerable computational complexity. Moreover, the solutions thus obtained in terms of recursion relations do not give much of structural insight in the solutions which, after all, are the noncommutative counterpart of something as simple and basic as plane waves. Why is this computation so difficult? It can be shown that just as in the commutative case the space of solutions of the massive Klein-Gordon Eq.~\eqref{eq:KleinGordon2} is generated, as representation of the quantum Poincar\'e algebra, by a rest state, which satisfies  \cite{Blohmann:2001b}
\begin{equation}
\label{eq:RestState1}  
  \partial^0 \tr \psi = \I m \psi
  \,,\qquad
  \partial^A \tr \psi = 0 \,,
\end{equation} 
where $A$ is a three vector index of spatial coordinates. In the commutative case we could now infer that $\psi$ is a function of the time coordinate $x_0$ alone, thus reducing Eq.~\eqref{eq:RestState1} to an ordinary differential equation in one coordinate $x_0$. Not so in the noncommutative case, where the partial derivatives of a function in $x_0$ depend on all coordinates $x^\mu$, which makes the solution of Eq.~\eqref{eq:RestState1} so involved. In analogy to partial differential equations in non-cartesian coordinates we will describe this mixing of dependencies by saying that within the noncommutative calculus the standard variables are not separated.

It is now obvious to ask, whether new variables can be found which separate the differential calculus, in such a way that Eq.~\eqref{eq:RestState1} becomes easily solvable. The main purpose of this paper is to show that this question has a positive answer.

\subsection{Structure, main results, notation}

The paper is organized as follows: In Sec.~\ref{sec:DiffCalc} the structure of the differential calculus on quantum Minkowski space is reviewed. We recall the commutation relations, general structure of the quantum Lorentz algebra, universal $\R$-matrix, and the definition of the coproduct of momenta in order to make this paper self contained and fix conventions and notation unambiguously.

Sec.~\ref{sec:Compute} contains the bulk of the paper: the computation of derivatives of arbitrary elements in the quantum Minkowski algebra. We first develop a notation and calculus in order to deal with functions of the algebra valued 4$\times$4-matrices which come from the action of the universal $L$-matrices appearing in the definition of the coproduct of momenta. The computation problem is here the calculation of powers of algebra valued matrices, which can then used to calculate the derivatives of arbitrary functions of the noncommutative coordinates. We do this calculation in two steps: In the first step we deal with polynomials of the non-central variables. The main result is in Proposition~\ref{th:SpaceDeriv}, Eqs.~\eqref{eq:SpaceDeriv} and can be written in the concise form of a chain rule, where the outer derivatives are partial Jackson derivatives. In the second step we calculate derivatives of functions of the central coordinates, time and four length. The partial derivatives of functions in the time generator are complicated and depend on all other generators. However, we are able to give a noncommutative coordinate transformation such that the derivatives of functions in the central generators are disentangled and take the form of a chain rule, too. This is the main result of the paper, presented in Proposition~\ref{th:CentralDeriv}, Eq.~\eqref{eq:ADeriv}. The results of these two steps can be combined to yield a compact formula for the derivatives of arbitrary elements of the algebra, given in Eq.~\eqref{eq:GeneralDeriv}.
 
In Sec.~\ref{sec:Wave} we apply the formulas for the derivations in order to the initial problem of finding the momentum eigenstate solutions of the massless and massive quantum Klein-Gordon equation, viewed as differential equations within the noncommutative differential calculus. We show that the separation of variables of the preceding section separates the differential equations completely. This reduces the partial differential equations to single variable equations involving Jackson derivatives, which are easily solved. The solutions are given in Eqs.~\eqref{eq:MasslessSolution} and \eqref{eq:MassiveSolution} in terms of $q$-exponential functions.
 
In Sec.~\ref{sec:Conclusion} we assess the results, give an outlook and references to further work in this direction of this paper. 

Appendix~\ref{sec:AppBasic} finally contains a few more technical results which are used in the paper. This makes this paper reasonably self-contained and will be useful for a reader who wishes to reproduce the calculations in detail.

Notation: Throughout the paper $q$ denotes a real deformation parameter, such that for $q \rightarrow 1$ we retrieve the undeformed, commutative limit. We will often use the abbreviations $\lambda := q - q^{-1}$ and symmetrized and standard quantum numbers
\begin{equation}
\label{eq:qNumbersDef}  
  [n] = \frac{q^n - q^{-n}}{q-q^{-1}}
  \,,\qquad
  \qNum{n} = \frac{q^{2n} - 1}{q^2-1} \,,
\end{equation}
where $n$ is a natural number. Repeated upper and lower indices are summed over $x_\mu x^\mu \equiv \sum_\mu x_\mu x^\mu$.

\section{Differential calculus on quantum Minkowski space}
\label{sec:DiffCalc}

\subsection{Quantum Minkowski space}

Quantum Minkowski space has first been constructed as braided product of two copies of the quantum plane \cite{Carow-Watamura:1990}. We will need the definition in terms of generators which are decomposed into the scalar time coordinate and the three vector of space coordinates.

\begin{Definition}
The complex algebra generated by the four generators $x_0$, $x_-$, $x_+$, and $x_3$, divided by the relations
\begin{equation}
\label{eq:XX-Rel1}
\begin{gathered}
  x_- x_0 = x_0 x_- \,,\quad x_+ x_0 = x_0 x_+ \,,\quad
  x_3 x_0 = x_0 x_3 \,,\\
  q^{-1} x_- x_3 - q x_3 x_- = -\lambda x_- x_0\,,\quad
  q^{-1} x_3 x_+ - q x_+ x_3 = -\lambda x_+ x_0 \,,\\ 
  x_- x_+ - x_+ x_- - \lambda x_3 x_3 = -\lambda x_3 x_0 \,,
\end{gathered}
\end{equation}
is called the algebra of quantum Minkowski space, denoted by $\Mink$.
\end{Definition}

Here, $x_0$ is the time coordinate, while $x_-$, $x_+$, and $x_3$ are the space coordinates with $\mathrm{su}_2$-weight indices. The commutation relations~\eqref{eq:XX-Rel1} can be written in more sophisticated forms, e.g., in terms of $q$-Clebsch-Gordon coefficients or in terms of an $R$-matrix as we will recall in the next section. 

The center of $\Xcal$ is generated by $x_0$ and the four-square
\begin{equation}
  \Xsquare := (x_0)^2  + q^{-1} x_- x_+ +  q x_+ x_- - (x_3)^2 \,.
\end{equation}
Since we will not use 2 as index, $\Xsquare$ cannot be confused with any one of the coordinates. The four-square can also be written as 
\begin{equation}
\label{eq:Xsquare} 
  \Xsquare = x_\mu x_\nu \eta^{\mu\nu}
\end{equation} 
in terms of the $q$-Minkowski metric defined by
\begin{equation}
  \eta^{00} = 1\,\quad \eta^{-+} = q^{-1} \,,\quad 
  \eta^{-+} = q\,\quad \eta^{33} = -1 \,,
\end{equation}
all other components vanishing. The metric defines upper four-vector indices by
\begin{equation}
  x^\mu := \eta^{\mu\nu} x_\nu \,.
\end{equation}
Note, that $\Xsquare = x_\mu x^\mu$ but $\Xsquare \neq x^\mu x_\mu$ because $\eta^{\mu\nu}$ is not symmetric.

\subsection{Quantum Lorentz algebra}

By construction $\Mink$ is a module algebra with respect to the quantum Lorentz algebra $\slC$. There are several essentially isomorphic forms of $\slC$ in use \cite{Schmidke:1991,Majid:1993,Lorek:1997a,Rohregger:1999}. In the chiral decomposition form $\slC$ is given, as algebra, by the tensor square of quantum $\Usu$,
\begin{equation}
\label{eq:qLorentzChiral}  
  \slC = \Uqsu \otimes \Uqsu \,.
\end{equation}
The definition of $\Uqsu$ is recalled in Appendix~\ref{sec:AppBasic}. The chiral decomposition implies that the irreducible representations of $\slC$ can be labelled by the highest weights of the two tensor factors. For example, the structure morphism of the four-vector representation is given on this chiral form by
\begin{equation}
  \rho^{(\frac{1}{2},\frac{1}{2})} 
  := \rho^{\frac{1}{2}} \otimes \rho^{\frac{1}{2}} \,,
\end{equation}
where $\rho^{\frac{1}{2}}$ is the fundamental representation of $\Uqsu$ with highest weight or spin $j=\frac{1}{2}$. 


The quantum Lorentz algebra is quasi-triangular. In the chiral form~\eqref{eq:qLorentzChiral} the universal $\R$-matrix of $\slC$ can be expressed in terms of the well known universal $\R$-matrix of $\Uqsu$. There are two inequivalent universal $\R$-matrices \cite{Majid},
\begin{equation}
\label{eq:Rdef}  
  \RI = \R^{-1}_{41}\R^{-1}_{31}\R_{24}\R_{23}\,, \qquad
  \RII = \R^{-1}_{41}\R_{13}\R_{24}\R_{23} \,,
\end{equation}
where $\R$ is the universal $\R$-matrix of $\Uqsu$ given in Eq.~\eqref{eq:Rmatrix}, and where we have used the tensor leg notation $\R_{23} = 1\otimes \R \otimes 1$ etc. The four-vector representations of the universal $\R$-matrices are denoted by
\begin{equation}
\begin{aligned}
  R_\mathrm{I}^{\mu\nu}{}_{\sigma\tau} 
  &:= \bigl( \rho^{(\frac{1}{2},\frac{1}{2})} 
  \otimes \rho^{(\frac{1}{2},\frac{1}{2})} \bigr)
  (\RI)^{\mu\nu}{}_{\sigma\tau} \,,\\
  R_\mathrm{II}^{\mu\nu}{}_{\sigma\tau} 
  &:= \bigl( \rho^{(\frac{1}{2},\frac{1}{2})} 
  \otimes \rho^{(\frac{1}{2},\frac{1}{2})} \bigr)
  (\RII)^{\mu\nu}{}_{\sigma\tau} \,.
\end{aligned}
\end{equation}
With the $R$-matrices the commutation relations~\eqref{eq:XX-Rel1} can be written in the compact form
\begin{equation}
\label{eq:XX-Rel2}  
  x_\sigma x_\tau =
  x_\mu x_\nu R_\mathrm{I}^{\mu\nu}{}_{\tau\sigma} \,.
\end{equation}
This implies that the $\Mink$-algebra has the Poincar\'e-Birkhoff-Witt property, because the $R$-matrices satisfy the Yang-Baxter equation. From the decomposition of the $R$-matrices into eigenspaces we can, furthermore, derive the relation
\begin{equation}
\label{eq:RR-Rel}  
  (1+ q \hat{R}_\mathrm{II})(1- \hat{R}_\mathrm{I}) = 0 \,,
\end{equation}
where the hat denotes swapping the two upper indices $\hat{R}^{\mu\nu}{}_{\tau\sigma} = R^{\nu\mu}{}_{\tau\sigma}$. We will need this formula later.

The disadvantage of the chiral form~\eqref{eq:qLorentzChiral} is that the inner tensor factors of the coproduct are twisted with the $\R$-matrix of $\Uqsu$. Since the $\R$-matrix is given by an infinite series, which exists only as formal power series, some of the generators of $\slC$ have a complicated coproduct which then exists only as formal power series, too. Therefore, it is also convenient to use the Drinfeld double form  \cite{Podles:1990}, where $\slC$ is represented as Drinfeld double of the Hopf-dually paired $\suq$ and $\SUq^\op$.

\subsection{Noncommutative differential calculus}

The generators of $q$-momenta are required to transform as four-vector with respect to $q$-Lorentz transformations. The fact that $q$-momenta transform in the same way as the coordinate generators of $q$-Minkowski space implies that they have to satisfy the same commutation relations as well. This leads to the definition of the $q$-momentum algebra as the algebra generated by $p_0$, $p_-$, $p_+$, $p_3$ with commutation relations~\eqref{eq:XX-Rel1}, replacing $x$ with $p$ everywhere.

Noncommutative partial derivatives 
\begin{equation}
\label{eq:Pact}
  \partial_\mu := \I p_\mu 
\end{equation} 
are defined by the action of the momentum generators on the $q$-Minkowski algebra of space functions. Partial derivatives of the coordinates have to be dimensionless numbers from the base field of complex numbers (including the deformation parameter $q$). Just as in the undeformed case the only tensor in the base field with two four-vector indices is the metric. We conclude that the action of partial derivatives on coordinates is given by
\begin{equation}
\label{eq:Pact2}  
  \partial_\mu \act x_\nu = \eta_{\mu\nu}
  \quad\Leftrightarrow\quad
  \partial^\mu \act x_\nu = \delta^\mu_\nu \,.  
\end{equation}
This action on the generators extends to the entire algebra of $q$-Minkowski space by the coproduct of the momenta \cite{Majid:1993},
\begin{equation}
\label{eq:Pcoprod}
\begin{split}
  \Delta(p^\mu) &:= p^\mu\otimes 1 +
  (\kappa\otimes 1)\, \RII^{-1} (1\otimes p^\mu) \RII \\
  &= p^\mu\otimes 1 + \kappa\, L^{\mu}{}_{\nu}
  \otimes p^{\nu} \,.
\end{split}
\end{equation}
Here $\RII$ is the second of the universal $\R$-matrices of the $q$-Lorentz algebra given in Eq.~\eqref{eq:Rdef}. The $L$-matrix is the four-vector half-representation of this $\R$-matrix
\begin{equation}
\label{eq:Ldef}
  L^{\mu}{}_{\nu} := \RII_{[1]}\,
    \rho^{(\frac{1}{2},\frac{1}{2})}(\RII_{[2]})^\mu{}_{\nu} \,,
\end{equation}
where we have used a Sweedler-like notation for $\RII = \RII_{[1]} \otimes \RII_{[2]}$. Finally, $\kappa$ is a group-like scaling operator, 
\begin{equation}
\label{eq:Scaling}  
  p_\mu \kappa = q\kappa p_\mu \,,\qquad
  x_\mu \kappa = q^{-1} \kappa x_\mu \,,\qquad
  \Delta(\kappa) = \kappa \otimes \kappa \,.
\end{equation}

The action~\eqref{eq:Pact2} of partial $q$-derivatives together with the coproduct~\eqref{eq:Pcoprod} defines a $q$-Lorentz-covariant differential calculus on $q$-Minkowski space. Letting $p^\mu$ act on $x_\nu f$, $f \in \Mink$, we obtain
\begin{equation}
  p^\mu \act (x_\nu f) = -\I \delta^{\mu}_\nu f 
  + q R_{\mathrm{II}}^{\nu'\mu}{}_{\nu\mu'} x_{\nu'} (p^{\mu'} \act f) \,,
\end{equation}
from which we can deduce the $q$-deformed Heisenberg commutation relations
\begin{equation}
  p^\mu x_\nu  - q R_{\mathrm{II}}^{\nu'\mu}{}_{\nu\mu'} x_{\nu'} p^{\mu'} 
  = -\I \delta^{\mu}_\nu  \,,
\end{equation}
which are of the type originally introduced by Wess and Zumino \cite{Wess:1991}.

In order to show that Eqs.~\eqref{eq:Pact2} and \eqref{eq:Pcoprod} define a differential calculus which is well defined on the algebra $\Mink$, that is, which does not depend on the ordering of the coordinate generators, we have to check if the action of the momenta is consistent with commutation relations~\eqref{eq:XX-Rel1} defining $\Mink$. This is best done using the $R$-matrix form~\eqref{eq:XX-Rel2},
\begin{equation}
\begin{split}
  p^\lambda \act (x_\sigma x_\tau -
  x_\mu x_\nu R_\mathrm{I}^{\mu\nu}{}_{\tau\sigma})
  &=
  (p^\lambda \act x_\nu x_\mu )( \delta^\nu_\tau \delta^\mu_\sigma -
  \hat{R}_\mathrm{I}^{\nu\mu}{}_{\tau\sigma}) \\
  &=
  x_\rho
  ( \delta^\lambda_\nu \delta^\rho_\mu +
  q \hat{R}_\mathrm{II}^{\lambda\rho}{}_{\nu\mu})
  ( \delta^\nu_\tau \delta^\mu_\sigma -
    \hat{R}_\mathrm{I}^{\nu\mu}{}_{\tau\sigma}) \\
  &= 0 \,,
\end{split} 
\end{equation} 
where in the last step we have used Eq.~\eqref{eq:RR-Rel}. 

In principle, the action~\eqref{eq:Pact2} of the partial derivatives on generators together with the coproduct~\eqref{eq:Pcoprod} of momenta defines the action of partial derivatives on arbitrary elements of quantum Minkowski space. Note, however, that this definition is recursive. It does not yield formulas for the partial derivatives of a basis of the noncommutative space algebra, such as the Poincar\'e-Birkhoff-Witt basis of ordered monomials,
\begin{equation}
  \mathcal{B}_{\mathrm{PBW}}
  := \{ (x_0)^{n_0} (x_-)^{n_-} (x_+)^{n_+} (x_3)^{n_3} 
  \,| \, n_0,n_-,n_+,n_3 \in \mathbb{N}_0 \} \,.
\end{equation}
To derive such a formula is one of the goals of the next section.

\section{Computation in noncommutative calculus}
\label{sec:Compute}

\subsection{\textit{L}-matrix calculus}

Let us first consider the simple case of a power of one single generator, $f = x_\alpha^n$, $\alpha \in \{0,-,+,3\}$, $n \in \mathbb{N}$. Using the coproduct~\eqref{eq:Pcoprod} we get by induction
\begin{equation}
\label{eq:MonDeriv1}
  \partial^\mu \act x_\alpha^n = \sum_{k=0}^{n-1}
  q^k (L^{\mu}{}_{\nu}\act x_\alpha^k)
  (\partial^{\nu} \act x_\alpha) x_\alpha^{n-k-1}  \,,
\end{equation}
where the factor $q^k$ comes from the action of the scaling operator~\eqref{eq:Scaling},
\begin{equation}
  \kappa \tr x_\alpha^k = \kappa_{(1)} x_\alpha^k S(\kappa_{(2)})
  = \kappa x_\alpha^k \kappa^{-1} = q^k x_\alpha^k \,.
\end{equation}
By construction, the $L$-matrix is a comultiplicative quantum matrix,
$\Delta(L^{\mu}{}_{\nu}) = L^{\mu}{}_{\sigma} \otimes
L^{\sigma}{}_{\nu}$. So we get for the action of $L$ on a power
\begin{equation}
\label{eq:Lalphapowers}  
  L^{\mu}{}_{\nu}\act x_\alpha^k
  = (L^{\mu}{}_{\sigma_1}\act x_\alpha )
    (L^{\sigma_1}{}_{\sigma_2}\act x_\alpha )
    \cdots
   (L^{ \sigma_{k-1} }{}_{\nu} \act x_\alpha )
   \,.
\end{equation}
Since each $L^{\sigma_k}{}_{\sigma_{k+1}}\act x_\alpha$ is an $\Xcal$-valued 4$\times$4-matrices, the right hand side is the matrix product of such matrices. This suggests to introduce an index-free notation:

\begin{Definition}
\label{th:indexfree}
Let $L^\mu{}_\nu$ be the $L$-matrix defined in Eq.~\eqref{eq:Ldef}, $\partial^\mu$ the partial derivatives defined in~\eqref{eq:Pact}, $f \in \Mink$ an element of the algebra of quantum Minkowski space. Then we denote by $L_f$ the $\Xcal$-valued 4$\times$4-matrix and by $\nabla f$ the $\Mink$-valued four-vector with entries
\begin{equation}
\label{eq:Lindices}
  (L_f)^{\mu}{}_{\nu} := L^{\mu}{}_{\nu} \act f
  \qquad\text{and}\qquad
  (\nabla f)^\mu := \partial^\mu \act f \,,
\end{equation}
respectively.
\end{Definition}

To illustrate this, we have for example
\begin{equation}
  (L^2_f \nabla f)^\mu = 
  (L^\mu{}_\nu \tr f) (L^\nu{}_\sigma \tr f) (\partial^\sigma \tr f) \,.
\end{equation} 
In this index-free notation Eq.~\eqref{eq:Lalphapowers} is written as $L \act x_\alpha^k = L_{x_\alpha}^k$, such that Eq.~\eqref{eq:MonDeriv1} becomes
\begin{equation}
\label{eq:MonDeriv2}
  \nabla x_\alpha^n 
  =\sum_{k=0}^{n-1} q^k L_{x_\alpha}^k
    (\nabla x_\alpha) x_\alpha^{n-k-1}  \,.
\end{equation}
We recall that by definition~\eqref{eq:Pact2} of the action of partial derivatives on the generators we have for the components of the gradient
\begin{equation}
  (\nabla x_\alpha)^\mu = \partial^\mu \act x_\alpha = \delta^\mu_\alpha \,,
\end{equation}
where $\mu, \alpha \in \{0,-,+,3\}$. 

The computationally difficult part of~\eqref{eq:MonDeriv2} is the evaluation of powers of the algebra-valued matrices $L_{x_\alpha}$. For example, for $(L_{x_0})^\mu{}_\nu$ we get with respect to four-vector indices $\mu,\nu \in \{0,-,+,3\}$
\begin{equation}
\label{eq:Lzero}
  L_{x_0} =
  \frac{1}{[2]}
  \begin{pmatrix}
    \frac{[4]}{[2]} x_0 & q\lambda x_- & q\lambda x_+ & q\lambda x_3 \\
    -\frac{\lambda}{q^2} x_+ & 
      \frac{[4]}{[2]} x_0 - \frac{\lambda^2}{q} x_3 & 0 & \lambda x_+\\
    \lambda x_- & 0 & 
      \frac{[4]}{[2]} x_0 - q\lambda^2 x_3 & -\lambda x_-\\
    \frac{\lambda}{q} x_3 & -q\lambda x_-& \frac{\lambda}{q} x_+ & 
      \frac{[4]}{[2]} x_0 - \lambda^2 x_3\\
  \end{pmatrix} .
\end{equation}
Little is known about the calculation of powers of general algebra valued matrices. Since the entries of this matrix do not commute we cannot resort to the usual methods of linear algebra. And in a brute force approach the number of terms in $L_{x_0}^n$ which have to be reordered to the Poin\-ca\-r\'e-Birk\-hoff-Witt form (or any other ordering) increases exponentially with $n$.

\subsection{Block decomposition of \textit{L}-matrices}

The computational problem of computing the powers of the $\Mink$-valued $L_{x_\alpha}$-matrices can be simplified by considering the $L$-matrix~\eqref{eq:Ldef} not with respect to four-vector indices running through $\{0,-,+,3\}$ as in Eq.~\eqref{eq:Lzero}, but with respect to pairs of indices of the two spin-$\frac{1}{2}$ representations of the chiral decomposition $\slC = \suq \otimes \suq$ of the $q$-Lorentz algebra, each labelled by the weights $\{-\frac{1}{2},
+\frac{1}{2}\} = \{-,+\}$.  We recall that the generators of quantum Minkowski space in this $(\frac{1}{2}, \frac{1}{2})$-spinor basis labelled by indices $\{--,-+,+-,++ \}$ are related to those in the four-vector basis by
\begin{equation}
\label{eq:Minkspinorbasis}
\begin{aligned}
  x_{--} &= [2]^{\frac{1}{2}} x_{-} \,,\\
  x_{-+} &= q^{\frac{1}{2}}(x_{3} - x_{0}) \,,\\
  x_{+-} &= q^{-\frac{1}{2}}x_{3} + q^{\frac{3}{2}}x_{0} \,,\\
  x_{++} &= [2]^{\frac{1}{2}} x_{+} \,.
\end{aligned}
\end{equation}
Using formula~\eqref{eq:Rdef} for the universal $\RII$-matrix in terms of the $\R$-matrix of $\Uqsu$ and definition~\eqref{eq:Ldef} of the $L$-matrix, we calculate the $L$-matrix with respect to the spinor basis,
\begin{equation}
\label{eq:Lboostrot}
\begin{split}
  (L)^{ij}{}_{kl} &= \bigl(\id\otimes\id\otimes
  \rho^{\frac{1}{2}}\otimes\rho^{\frac{1}{2}} \bigr)
  (\RII)^{ij}{}_{kl} \\
  &= \bigl( L^{\frac{1}{2}}_- \bigr)^{j}{}_{j'} \bigl(
  L^{\frac{1}{2}}_+ \bigr)^{i}{}_{i'} \otimes \bigl( L^{\frac{1}{2}}_+
  \bigr)^{j'}{}_{l}
  \bigl( L^{\frac{1}{2}}_+ \bigr)^{i'}{}_{k}\\
  &= B^{j}{}_{l} (L^{\frac{1}{2}}_+)^{i}{}_{k} \,,
\end{split}
\end{equation}
where all indices run through $\{-\frac{1}{2}, +\frac{1}{2}\} = \{-,+\}$.  Here $B=(B^{j}{}_{l})$ is the matrix of ``boosts'' which generate the $\SUq^\op$ Hopf subalgebra of the Drinfeld double form of the $q$-Lorentz algebra
\cite{Podles:1990}, and $L^{1/2}_+$ the $L_{+}$-matrix
\cite{Schmuedgen} of the $\suq$ subalgebra of rotations,
\begin{equation}
  B  = \begin{pmatrix} a & b \\ c & d \end{pmatrix} \,,\quad
  L_{+}^{\frac{1}{2}} = \begin{pmatrix}
    K^{-\frac{1}{2}} &
    q^{-\frac{1}{2}}\lambda K^{-\frac{1}{2}}E \\
              0      & K^{\frac{1}{2}}   \end{pmatrix},
\end{equation}
with respect to the $\{-,+\}$ basis.  The explicit form of the boosts
in the chiral decomposition of the $q$-Lorentz algebra is given in Eq.~\eqref{eq:boostdef}. 

The decomposition~\eqref{eq:Lboostrot} of the $L$-matrix into boosts
and rotations yields a natural block matrix decomposition of the $\Mink$-valued 4$\times$4-matrices $L_{x_\alpha}$ as defined in Eq.~\eqref{eq:Lindices}. For the four-vector index $\alpha$ it turns out to be convenient to use instead of $x_3$ the light cone coordinate 
\begin{equation}
  x_\Lind := x_3 -x_0 \,,
\end{equation}
which is essentially the $x_{-+}$-coordinate in the spinor basis~\eqref{eq:Minkspinorbasis} of quantum Minkowski space. And for the block decomposition it is convenient to define $\Xcal$-valued 2$\times$2-matrices by the action of the boosts
\begin{equation}
  (B_f)^i{}_{j} := B^i{}_{j}\tr f \,, \qquad f\in \Xcal \,,
\end{equation}
using the same index-free notation as introduced in Eq.~\eqref{eq:Lindices}. Since $B^i{}_{j}$ is a comultiplicative quantum matrix, $\Delta(B^i{}_{j}) = B^i{}_{k} \otimes B^k{}_{j}$, the matrices $B_{x_\alpha}$ satisfy the same commutation relations~\eqref{eq:XX-Rel1} as $x_\alpha$. Explicitly, the $B_{x_\alpha}$-matrices can be computed to
\begin{equation}
\label{eq:AppBoostAct}
\begin{aligned}
  \Mat{B}_{x_0} &=
    \begin{pmatrix}
      \frac{[4]}{[2]^2}x_0 + \frac{\lambda}{q[2]} x_3  &
      q^{-\frac{1}{2}}\lambda [2]^{-\frac{1}{2}} x_+ \\
      -q^{\frac{1}{2}}\lambda [2]^{-\frac{1}{2}} x_- &
      \frac{[4]}{[2]^2} x_0 - \frac{q \lambda}{[2]} x_3
    \end{pmatrix} \\
  \Mat{B}_{x_-} &=
    \begin{pmatrix} x_- \quad &
      q^{-\frac{1}{2}}\lambda [2]^{-\frac{1}{2}} x_\Lind \\
      0 & x_-
    \end{pmatrix} \\
  \Mat{B}_{x_+} &=
    \begin{pmatrix} x_+ & 0 \\
      -q^{\frac{1}{2}}\lambda [2]^{-\frac{1}{2}} x_\Lind & \quad x_+
    \end{pmatrix}\\
  \Mat{B}_{x_\Lind} &=
    \begin{pmatrix} q^{-1}x_\Lind & 0 \\
    0 & qx_\Lind \end{pmatrix},
\end{aligned}
\end{equation}
with respect to spin-$\frac{1}{2}$ indices running through $\{-,+\}$. In terms of these matrices the $L_{x_\alpha}$-matrices can be expressed as
\begin{equation}
\label{eq:Lblock}
\begin{aligned}
  \Mat{L}_{x_0}
  &= \begin{pmatrix} \Mat{B}_{x_0} & 0 \\
    0 & \Mat{B}_{x_0} \end{pmatrix} \\
  \Mat{L}_{x_-}
  &= \begin{pmatrix} q\Mat{B}_{x_-} &
    q^{-\frac{1}{2}}\lambda [2]^{\frac{1}{2}} \Mat{B}_{x_3} \\
    0 & q^{-1}\Mat{B}_{x_-} \end{pmatrix} \\
  \Mat{L}_{x_+}
  &= \begin{pmatrix} q^{-1}\Mat{B}_{x_+} & 0 \\
    0 & q \Mat{B}_{x_+} \end{pmatrix} \\
  \Mat{L}_{x_\Lind}
  &= \begin{pmatrix} \Mat{B}_{x_\Lind} &
    q^{-\frac{1}{2}} \lambda [2]^{\frac{1}{2}} \Mat{B}_{x_+} \\
    0 & \Mat{B}_{x_\Lind} \end{pmatrix},
\end{aligned}
\end{equation}
with respect to the index structure $\{--,-+,+-,++\}$.

\subsection{Calculating powers of \textit{L}-matrices}

In the form of Eq.~\eqref{eq:Lblock} the $L_{x_\alpha}$ matrices are upper block triangular. This makes the computation of their powers a lot easier. $\Mat{L}_{x_0}$ and $\Mat{L}_{x_+}$ are even block diagonal, so we immediately get
\begin{equation}
\label{eq:L+Powers}  
  \Mat{L}_{x_0}^n
  = \begin{pmatrix} \Mat{B}_{x_0}^n & 0 \\
    0 & \Mat{B}_{x_0}^n \end{pmatrix} \,, \quad
  \Mat{L}_{x_+}^n
  = \begin{pmatrix} q^{-n}\Mat{B}_{x_+}^n & 0 \\
    0 & q^{n} \Mat{B}_{x_+}^n \end{pmatrix}. 
\end{equation}
For the calculation of $\Mat{L}_{x_-}^n$ and $\Mat{L}_{x_\Lind}^n$ we have to use the commutation relations~\eqref{eq:XX-Rel1} for the $B_{x_\alpha}$ matrices. We thus get by induction
\begin{align}
\label{eq:L-LindPowers}  
  \Mat{L}_{x_-}^n
  &= \begin{pmatrix} q^{n} \Mat{B}_{x_-}^n &
    \; q^{-\frac{1}{2}}\lambda [2]^{\frac{1}{2}}
    \left( q^{n-1}\frac{[2n]}{[2]} \Mat{B}_{x_\Lind} + [n]\Mat{B}_{x_0}
      \right) \Mat{B}_{x_-}^{n-1}  \\
    0 & q^{-n} \Mat{B}_{x_-}^n \end{pmatrix} \notag \\
  \Mat{L}_\Lind^n 
  &= \begin{pmatrix} \Mat{B}_{x_\Lind}^n &
    \; q^{-\frac{1}{2}} \lambda [2]^{\frac{1}{2}}
    q^{n-1} [n] \Mat{B}_{x_+} \Mat{B}_{x_\Lind}^{n-1}  \\
    0 & \Mat{B}_{x_\Lind}^n \end{pmatrix}.
\end{align}
Now we have expressed the powers of $L_{x_\alpha}$ matrices in terms of powers of the $B_{x_\alpha}$ matrices. It remains to calculate the powers of the $B_{x_\alpha}$-matrices. Since $B_{x_\Lind}$ is diagonal we immediately get
\begin{equation}
\label{eq:BLindPowers}  
  \Mat{B}_{x_\Lind}^n
  = \begin{pmatrix} q^{-n}x_\Lind^n & 0 \\
    0 & q^n x_\Lind^n \end{pmatrix} \,. 
\end{equation}
For the powers of the matrices $B_{x_\pm}$, which are block triangular, we use again the commutation relations~\eqref{eq:XX-Rel1} to obtain by induction
\begin{equation}
\label{eq:B-+Powers}
\begin{aligned}
  \Mat{B}_{x_-}^n
  &= \begin{pmatrix} x_-^n &
    \; q^{-\frac{1}{2}}\lambda [2]^{-\frac{1}{2}}
    q^{n-1}[n] x_\Lind x_-^{n-1}  \\
    0 & x_-^n \end{pmatrix} \\
  \Mat{B}_{x_+}^n
  &= \begin{pmatrix} x_+^n & 0  \\
    -q^{\frac{1}{2}}\lambda [2]^{-\frac{1}{2}}
    q^{1-n}[n] x_\Lind x_+^{n-1} & x_+^n \end{pmatrix}.
\end{aligned}
\end{equation}
The calculation of the powers of $B_{x_0}$ is more difficult because it is not triangular. Using $q$-Pauli matrices~\eqref{eq:sigmaDef}, $B_{x_0}$ of Eqs.~\eqref{eq:AppBoostAct} can be written more compactly as
\begin{equation}
\label{eq:B0compact}
  \Mat{B}_{x_0} := \frac{[4]}{[2]^2} x_0 -
  \frac{\lambda}{[2]} \tilde{\sigma}\!_{A} x^A \,,
\end{equation}
where the index $A$ is summed over $\{-,3,+\}$. The $q$-Pauli matrices satisfy the relation
\begin{equation}
\label{eq:sigmatilderel}
  \tilde{\sigma}\!_A \tilde{\sigma}_B
  = g_{BA} - \tilde{\sigma}_C\, \varepsilon_B{}^C{}_A \,,
\end{equation}
where the quantum metric and epsilon tensor are defined by quantum Clebsch-Gordon coefficients as 
\begin{equation}
\begin{aligned}
  g^{AB} 
  &:= -\sqrt{[3]} \,\, C_q(1,1,0;A,B,0) \\ 
  \varepsilon^{AB}{}_{C} 
  &:= -\sqrt{\frac{[4]}{[2]}} \,\, C_q(1,1,1;A,B,C) \,.
\end{aligned}
\end{equation}
Contracting Eq.~\eqref{eq:sigmatilderel} with $x^A x^B$ and using the commutation relations~\eqref{eq:XX-Rel1} of the coordinates, which can be written as $x_A x_B \varepsilon^{AB}{}_C = -\lambda x_C x_0$, we
derive
\begin{equation}
\begin{split}
  \tilde{\sigma}\!_{A} \tilde{\sigma}\!_{B} x^A x^B 
  &= (g_{BA} - \tilde{\sigma}_C\, \varepsilon_B{}^C{}_A) x^A x^B \\
  &= x_A x^A + \lambda x_0\,\tilde{\sigma}_A\, x^A \,.
\end{split}
\end{equation}
This can be used to compute the square of Eq.~\eqref{eq:B0compact},
\begin{equation}
\begin{split}
\label{eq:Msquare}
  \Mat{B}_{x_0}^2
  &= \frac{[4]^2}{[2]^2} (x_0)^2 - 2\frac{\lambda[4]}{[2]^3} x_0\,
    \tilde{\sigma}\!_{A} x^A + \frac{\lambda^2}{[2]^2}\,
    \tilde{\sigma}\!_{A} \tilde{\sigma}\!_{B} x^A x^B \\
  &= \frac{[4]^2}{[2]^4} (x_0)^2 - 2\frac{\lambda[4]}{[2]^3} x_0\,
    \tilde{\sigma}\!_{A} x^A 
  + \frac{\lambda^2}{[2]^2}\,
    (x_A x^A + \lambda x_0\,\tilde{\sigma}_A\, x^A)\\
  &= [2]x_0 \Mat{B}_{x_0} - (x_0)^2 
     - \frac{\lambda^2}{[2]^2} \Xsquare  \,.
\end{split}
\end{equation}
Together with Eq.~\eqref{eq:L+Powers} we conclude that $L_{x_0}$ satisfies the polynomial equation
\begin{equation}
\label{eq:L0char}
\begin{split}
  0 
  &= L_{x_0}^2 - [2]x_0 L_{x_0} 
  + \Bigl( (x_0)^2 + \frac{\lambda^2}{[2]^2} \Xsquare \Bigr) \\
  &= L_{x_0}^2 - 2 b L_{x_0} + c \,,
\end{split}
\end{equation}
where we have introduced the abbreviations
\begin{equation}
\label{eq:ABDef}
  b := \frac{[2]}{2} x_0
  \,,\qquad
  c:= (x_0)^2 + \frac{\lambda^2}{[2]^2} \Xsquare \,,
\end{equation}
which are both in the center of $\Mink$. Eq.~\eqref{eq:L0char} enables us to derive formulas for arbitrary powers of $L_{x_0}$. This is most elegantly done by considering the generating function of the powers $(1 - L_{x_0} z )^{-1} = 1 + L_{x_0} z + L_{x_0}^2 z^2 + \ldots$, where $z$ is a formal parameter. First we rewrite~\eqref{eq:L0char} as
\begin{equation}  
  (1- L_{x_0}^2 z^2) - 2 b z (1 - L_{x_0} z) - (1 - 2b z  + c z^2) = 0
\end{equation}
and divide it by $(1 - L_{x_0} z)$ and $(1 - 2b z  + c z^2)$, yielding
\begin{equation}  
  \frac{1}{1- L_{x_0} z}
  = \frac{1 + (L_{x_0} -2 b) z }{1 - 2 b z + c z^2 } \,.
\end{equation} 
The right hand side can be expanded in powers of $z$, observing that it contains the generating function for Chebyshev polynomials of the second kind, 
\begin{equation}
  \frac{1}{1 - 2y t + t^2} = \sum_{n=0}^\infty U_n(y)t^n \,,
\end{equation} 
where we have to set $y = b/\sqrt{c}$ and $t= \sqrt{c}\, z$. We thus obtain
\begin{equation}
\label{eq:L0Powers}
\begin{split}
  L_{x_0}^n 
  &=
    U_{n}\bigl( b/\sqrt{c} \bigr)\, c^{\frac{n}{2}}
  + U_{n-1} \bigl( b / \sqrt{c} \bigr)\, c^{\frac{n-1}{2}} (L_{x_0} -2b)  \\
  &= 
  - U_{n-2}\bigl( b/\sqrt{c} \bigr)\, c^{\frac{n}{2}} 
  + U_{n-1} \bigl( b / \sqrt{c} \bigr)\, c^{\frac{n-1}{2}} L_{x_0} 
\end{split}
\end{equation}
for $n \geq 2$. The fact that $U_n(y)$ is a polynomial of degree $n$ which contains either only even or only odd powers of $y$ implies that $U_n(1/y)y^n$ is a quadratic polynomial in $y$. Therefore, the right hand side of Eq.~\eqref{eq:L0Powers} contains only even positive powers of $\sqrt{c}$, that is, does not contain inverses or square roots of $c$ and, hence, yields well-defined elements of the algebra $\Xcal$. Together with Eqs.~\eqref{eq:L+Powers}, 
\eqref{eq:L-LindPowers}, \eqref{eq:BLindPowers}, and \eqref{eq:B-+Powers} this completes our calculation of powers of the $L_{x_\alpha}$ matrices.

\subsection{Calculating derivatives}

\subsubsection{Choice of basis and partial Jackson derivatives}

With the formulas for the powers of the $L_\alpha$ matrices the derivatives of arbitrary elements $f \in \Xcal$ could be computed by calculating the derivatives of the Poincar\'e-Birkhoff-Witt basis of $\Mink$. However, it turns out to be more convenient to work in a slightly different basis. Observe that from definition~\eqref{eq:Xsquare} of the $q$-Lorentz invariant length and the commutation relations~\eqref{eq:XX-Rel1} we can deduce the relation 
\begin{equation}
  [2] x_- x_+ = 
  \Xsquare  + q^2 (x_\Lind)^2 + q \lambda x_0 x_\Lind \,.
\end{equation} 
Using this relation, all products of equal powers of $x_-$ and $x_+$ in the Poincar\'e-Birkhoff-Witt basis can be substituted by powers of $\Xsquare$. We thus obtain a basis of ordered monomials which contains powers of $\Xsquare$ and either powers of $x_-$ or $x_+$. To be more precise, we have a basis of $\Mink$,
\begin{equation}
  \mathcal{B} = \mathcal{B}_- \cup \mathcal{B}_{+} \,,
\end{equation}
where
\begin{equation}
\label{eq:GeneralMon}
\begin{aligned}
  \mathcal{B}_-
  &:= \{ (x_\mu x^\mu)^i (x_0)^{j} (x_\Lind)^{k} (x_-)^{l}
  \,| \, i,j,k,l \in \mathbb{N}_0 \} \\
  \mathcal{B}_{+}
  &:= \{ (x_\mu x^\mu)^i (x_+)^{l} (x_0)^{j} (x_\Lind)^{k} 
  \,| \, i,j,k,l \in \mathbb{N}_0 \} \,.  
\end{aligned}
\end{equation}
Note, that $\mathcal{B}_-$ and $\mathcal{B}_+$ have a non-empty intersection, the monomials for which $l=0$. But this will not be a problem here. Let us introduce a more suggestive notion for sums of monomials in a particular order:
\begin{Definition}
\label{th:Order}
  Let $x_1,\ldots,x_n \in \Xcal$ be elements of the algebra of quantum Minkowski space. Then we will denote by $f(x_1,\ldots,x_n)$ linear combinations of ordered products of monomials in $x_1,\ldots,x_n \in \Xcal$. That is, 
  \begin{equation}
    f(x_1,\ldots,x_n) \in 
    \mathrm{Span}\{  (x_1)^{k_1} \cdots (x_n)^{k_n} 
    \,| \, k_1,\ldots,k_n \in \mathbb{N}_0 \} \,.
  \end{equation} 
\end{Definition} 
With this notation the statement that $\mathcal{B}_- \cup \mathcal{B}_+$ is a basis of quantum Minkowski space can be written as
\begin{equation}
  f(x_0,x_+,x_\Lind,x_-)
  =
  f_-(\Xsquare, x_0, x_\Lind, x_-) +
  f_+(\Xsquare, x_0, x_+, x_\Lind) \,, 
\end{equation}
decomposing an arbitrary element of $\Mink$ into two parts containing powers in either $x_-$ or $x_+$.

It is important to have a notation which keeps track of the ordering of generators because the formulas for derivatives are most elegantly expressed in terms of partial Jackson derivatives which depend on the ordering. Recall, that the Jackson derivative or $q$-derivative of a function $f = f(x_\alpha)$ in a single variable $x_\alpha$ is defined as difference quotient,
\begin{equation}
  \frac{\partial_{q^2} f}{\partial_{q^2}x_\alpha}
  := \frac{f(q^2 x_\alpha) - f(x_\alpha)}{q^2 x_\alpha - x_\alpha} \,.
\end{equation}
For monomials $f(x_\alpha) = (x_\alpha)^k$, this yields
\begin{equation}
  \frac{\partial_{q^2} (x_\alpha)^k}{\partial_{q^2}x_\alpha}
  = \qNum{k} (x_\alpha)^{k-1} \,,
\end{equation}
where $\qNum{k}$ is the usual $q$-number defined in Eq.~\eqref{eq:qNumbersDef}. This formula naturally generalizes to a partial Jackson derivative on monomials of several noncommutative variables. 

\begin{Definition}
\label{th:Jackson}
Let $x_1,\ldots,x_n \in \Mink$ be elements of the algebra of quantum Minkowski space (or any other algebra). We define the partial Jackson derivative with respect to one of these elements $x_\alpha$ on ordered monomials by
\begin{equation}
\label{eq:partialJackson}
  \frac{\partial_{q^2}}{\partial_{q^2}x_\alpha}
  \bigl( (x_1)^{k_1} \cdots (x_n)^{k_n} \bigr)
  := \qNum{k_\alpha}
    (x_1)^{k_1} \cdots (x_\alpha)^{k_\alpha-1} \cdots (x_n)^{k_n} \,,
\end{equation}
and extend it linearly to arbitrary linear combinations of such monomials. This defines partial Jackson derivatives on general functions $f = f(x_1,\ldots,x_n)$ in the sense of Definition~\ref{th:Order}.
\end{Definition}

It must be emphasized that the partial Jackson derivatives depend on the particular order and are not well-defined on the algebra. For example, within the algebra we have the commutation relation $x_\Lind x_+ - q^2 x_+ x_\Lind = 0$, but
\begin{equation}
  \frac{\partial_{q^2}}{\partial_{q^2}x_+}
  (x_\Lind x_+ - q^2 x_+ x_\Lind ) = x_\Lind - q^2 x_\Lind 
  \neq 0 \,.
\end{equation} 
This is why we have defined the notation $f = f(x_+, x_\Lind)$ to not just denote elements of the algebra but to denote in addition a specific ordering of the variables. With this notation the action of partial Jackson derivatives is defined unambiguously.

We break up the calculation of derivatives in two steps. First, we calculate the derivatives of non-central variables, $x_-$, $x_\Lind$, $x_+$. Working within the basis $\mathcal{B}$ we can limit the considerations in this first step to polynomial functions  $f = f(x_\Lind,x_- )$ and $f = f(x_+,x_\Lind)$ which depend either on $x_-$ or on $x_+$. This calculation is rather straightforward. In a second step we consider polynomial functions $f = f(\Xsquare, x_0)$ of the center of $\Xcal$. This calculation is much more involved.

\subsubsection{Functions of non-central coordinates}

We begin the calculation of the derivatives of the basis $\mathcal{B}$ by reading off Eqs.~\eqref{eq:AppBoostAct} and \eqref{eq:Lblock} that the $L_\alpha$ matrices possess some algebra valued eigenvalues, 
\begin{equation}
\label{eq:Eigen}
\begin{aligned}
  L_A \nabla x_A &= q x_A \nabla x_A \\
  L_\Lind \nabla x_- &= q^{-1} x_\Lind \nabla x_- \\
  L_+ \nabla x_\Lind &= q^{-1} x_+ \nabla x_\Lind \,,
\end{aligned}
\end{equation} 
for $A \in \{-, +, \Lind \}$ (no summation over $A$). We recall that $\nabla x_\alpha$ is the index free notation for $(\nabla x_\alpha)^\mu = \partial^\mu \tr x_\alpha = \delta^\mu_\alpha$.

Inserting the first of these eigenvalue equations into Eq.~\eqref{eq:MonDeriv2} we obtain for powers of the generators
\begin{equation}
\label{eq:Der1}
  \nabla x_A^n = \sum_{k=0}^n
  q^{2k} x_A^k (\nabla x_\alpha) x_\alpha^{n-k-1} 
  = \qNum{n} \, x_A^{n-1} \nabla x_\alpha \,,
\end{equation}
for $A \in \{-, +, \Lind \}$, where $\qNum{n}$ denotes the quantum number~\eqref{eq:qNumbersDef}. From Eq.~\eqref{eq:Der1} we can deduce that the derivative of any polynomial $f=f(x_A)$ in a single one of the generators $x_A$ can be expressed in terms of the Jackson derivative by 
\begin{equation}
\label{eq:Der2}
  \nabla f(x_A)   
  = \frac{\partial_{q^2} f }{\partial_{q^2}x_A} \,\nabla x_A 
  \quad\Leftrightarrow\quad
  \partial^\mu \tr f(x_A)   
  = \frac{\partial_{q^2} f }{\partial_{q^2}x_A} \, \delta_A^\mu 
  \,.
\end{equation}
Using the second of Eqs.~\eqref{eq:Eigen} we get
\begin{equation}
\label{eq:Der3}
\begin{split}
  \nabla (x_\Lind^k x_-^n )
  &= (\nabla x_\Lind^k) x_-^n + q^k L_\Lind^k (\nabla x_-^n) \\
  &= \qNum{k} x_\Lind^{k-1}  x_-^n \nabla x_\Lind 
  +  x_\Lind^k \,\qNum{n} x_-^{n-1} \nabla x_- \,,
\end{split}
\end{equation}
and analogously for the derivative of $x_+^k x_\Lind^n$,
\begin{equation}
\label{eq:Der4}
  \nabla (x_+^k x_\Lind^n )
  = \qNum{k} x_+^{k-1}  x_\Lind^n \nabla x_+ 
  + x_+^k \,\qNum{n} x_\Lind^{n-1} \nabla x_\Lind \,.
\end{equation}
From Eqs.~\eqref{eq:Der3} and \eqref{eq:Der4} we can deduce the following result:

\begin{Proposition}
\label{th:SpaceDeriv}
Let $f=f(x_\Lind, x_-)$ and $f=f(x_+,x_\Lind)$ be ordered polynomials in the notation of Definition~\ref{th:Order}. Then their partial derivatives are expressed in terms of partial Jackson derivatives as
\begin{subequations}
\label{eq:SpaceDeriv}
\begin{align}
  \nabla f(x_\Lind, x_-)   
  &= \frac{\partial_{q^2} f }{\partial_{q^2}x_\Lind} \nabla x_\Lind
   + \frac{\partial_{q^2} f }{\partial_{q^2}x_-} \nabla x_- \\
  \nabla f(x_+, x_\Lind )
  &= \frac{\partial_{q^2} f }{\partial_{q^2}x_+} \nabla x_+
   + \frac{\partial_{q^2} f }{\partial_{q^2}x_\Lind} \nabla x_\Lind
  \,.
\end{align}
\end{subequations}
\end{Proposition}

\subsubsection{Functions of four-length}

Next, we calculate the derivative of powers of the coordinate four-square $\Xsquare$. Since $\Xsquare$ is a Lorentz scalar the $L$-matrix acts with the antipode $\varepsilon$,
\begin{equation}
  L^\mu{}_\nu \tr \Xsquare
  = \varepsilon(L^\mu{}_\nu) \, \Xsquare
  = \delta^\mu_\nu \, \Xsquare \,.
\end{equation}
By the same reasoning which led to Eq.~\eqref{eq:MonDeriv2} we obtain for powers of the four-square
\begin{equation}
\label{eq:Der5}  
\begin{split}  
  \nabla(\Xsquare)^n
  &= \sum_{k=0}^{n-1} q^{2k} (\Xsquare)^{n-1} \nabla\Xsquare \\
  &= \qNum{n} (\Xsquare)^{n-1} \nabla\Xsquare
  \,.
\end{split}
\end{equation}
Again, for a general function $f = f(\Xsquare)$ this can be written in terms of a Jackson derivative by
\begin{equation}
\label{eq:XsquareDer}  
  \nabla f(\Xsquare)   
  = \frac{\partial_{q^2} f }{\partial_{q^2}(\Xsquare)}
  \nabla\Xsquare \,.
\end{equation}
It remains to compute the derivative of $\Xsquare$, which for representation theoretic reasons we expect to be proportional to the coordinate four vector $x^\mu$. After lengthy calculations we, indeed, find
\begin{equation}
\label{eq:SquareDeriv}
  (\nabla \Xsquare)^\mu = q^{-1} [2] x^\mu \,.
\end{equation}

\subsubsection{Functions of time}

Like the calculation of powers of $L_{x_0}$, the calculation of derivatives of powers of the time coordinate is more difficult. In order to calculate the derivative of $(x_0)^n$ we first observe that, since $x_0$ is central in the space algebra we can write Eq.~\eqref{eq:MonDeriv2} in the form
\begin{equation}
\label{eq:Lquotient}  
  \nabla x_0^n 
  = \frac{(q L_{x_0})^n - (x_0)^n}{qL_{x_0} - x_0} \nabla x_0 \,,
\end{equation}
which is similar to a $q$-difference quotient. We can get rid of the matrix in the denominator,
\begin{equation}
\begin{split}
  \frac{(q L_{x_0})^n - (x_0)^n}{qL_{x_0} - x_0}
  &= \frac{(q^n L_{x_0}^n  - x_0^n)(q^{-1}L_{x_0} - x_0)}
          {(qL_{x_0} - x_0)(q^{-1}L_{x_0} - x_0)} \\
  &=-\frac{[2](q^{n-1}L_{x_0}^{n+1}  - q^n L_{x_0}^n x_0 
           - q^{-1}L_{x_0} x_0^n - x_0^{n+1})}
          {\lambda^2 \Xsquare} \,,
\end{split} 
\end{equation}
where in the second step we have used Eq.~\eqref{eq:L0char}. The powers of $L_{x_0}$ in the numerator have been calculated in Eq.~\eqref{eq:L0Powers}. From explicit expression~\eqref{eq:Lzero} of $L_{x_0}$ we deduce
\begin{equation}
  (L_{x_0} \nabla x_0)^\mu = q x_0\, \delta^\mu_0 
  - \frac{\lambda}{q[2]} x^\mu \,.
\end{equation}
Putting things together we obtain
\begin{equation}
\label{eq:TimeDeriv}
\begin{split}
  (\nabla x_0^n)^\mu
  &= \frac{[2]}{\lambda}
  \biggl( 
  \frac{U_n c^{\frac{n}{2}} q^{n-2} - U_{n-1}c^{\frac{n-1}{2}} q^{n-1}  x_0 
  - x_0^n}{\Xsquare} 
  \biggr) x^\mu \\
  &\quad +
  U_{n-1}c^{\frac{n-1}{2}} q^{n-1} \delta^\mu_0 \,,
\end{split} 
\end{equation}
where $U_n = U_n(b/\sqrt{c})$ denotes the Chebyshev polynomials of the second kind with the same argument as in Eq.~\eqref{eq:L0Powers}.

In order to generalize the formulas to general functions $f = f(x_0)$ we first note that we can deduce from Eq.~\eqref{eq:Lquotient}
\begin{equation}
\label{eq:Lquotient2}  
  \nabla f(x_0) 
  = \frac{f(q L_{x_0}) - f(x_0)}{qL_{x_0} - x_0} \nabla x_0 \,,
\end{equation}
which could be viewed as matrix valued generalization of the Jackson derivative. In order to compute this expression we need to generalize formula~\eqref{eq:L0Powers} for the powers of $L_{x_0}$ to general functions. The result is given by:

\begin{Proposition}
\label{th:functionL}
  Let $L_{x_0}$ be the $\Mink$-valued 4$\times$4-matrix given by Eq.~\eqref{eq:Lzero}. Let $b$ and $c$ be defined as in Eq.~\eqref{eq:ABDef}. Define
\begin{equation}
\label{eq:alphaDef}  
  \tau_{\pm} := b \pm \sqrt{b^2-c} 
  = \frac{1}{2} \left( [2] x_0 
  \pm \lambda \sqrt{(x_0)^2 - \tfrac{4}{[2]^2} (\Xsquare)} \right) \,,
\end{equation}
  and
\begin{equation}
\label{eq:PiDef}
  \Pi_\pm := 
  \frac{1}{2}\! \left( 1 \pm \frac{L_{x_0} - b}{\sqrt{b^2 - c}} 
  \right) \,. 
\end{equation}
Then for any polynomial function $f$ in one variable we have
\begin{equation}
\label{eq:fLzero}
  f(L_{x_0}) = f(\tau_+) \Pi_+ + f(\tau_-) \Pi_-\,.
\end{equation} 
\end{Proposition}

\begin{proof}
First, we note that Chebyshev polynomials of the second kind can be written as
\begin{equation}
  U_{n-1}(x) = \frac{(x + \sqrt{x^2-1})^{n}- (x - \sqrt{x^2-1})^{n}}
  {2\sqrt{x^2-1}} \,.
\end{equation} 
For $x = b/\sqrt{c}$ this identity takes the form
\begin{equation}
\label{eq:Urel1}  
  U_{n-1}(b/\sqrt{c}) c^{\frac{n-1}{2}}
  = \frac{\tau_+^n - \tau_-^n}{\tau_+ - \tau_-} \,,
\end{equation} 
with $\tau_\pm$ defined as in the proposition. Inserting \eqref{eq:Urel1} into \eqref{eq:L0Powers} yields for an arbitrary function
\begin{equation}
\label{eq:L0function}
  f(L_{x_0}) 
  = \frac{f(\tau_+) - f(\tau_-)}{\tau_+ - \tau_-} (L_{x_0} - b)
  + \tfrac{1}{2}\bigl( f(\tau_+) + f(\tau_-) \bigr) \,,
\end{equation} 
which can be written using $\Pi_\pm$ in the form of Eq.~\eqref{eq:fLzero}.
\end{proof}

The $\Mink$-valued matrices $\Pi_\pm$ arise here naturally because they are orthogonal, complementary idempotents,
\begin{equation}
  \Pi_\pm^2 = \Pi_\pm \,,\quad \Pi_+ \Pi_- = 0 =  \Pi_- \Pi_+ 
  \,,\quad \Pi_+ + \Pi_- = 1 \,.
\end{equation}
This property ensures that Eq.~\eqref{eq:fLzero} is consistent with the algebra structure of functions, because it can be immediately verified that $(fg)(L_{x_0}) = f(L_{x_0}) g(L_{x_0})$.

We can now insert Eq.~\eqref{eq:L0function} into Eq.~\eqref{eq:Lquotient2} to obtain
\begin{equation}
\label{eq:Lquotient3}  
  \nabla f(x_0) = 
  \frac{f(q \tau_+) - f(x_0)}{q\tau_+ - x_0} \Pi_+ \nabla x_0
 +\frac{f(q \tau_-) - f(x_0)}{q\tau_- - x_0} \Pi_- \nabla x_0 
\end{equation} 
for the derivative of an arbitrary function of the time coordinate. For this formula to be explicit it remains to calculate
\begin{equation}
\label{eq:PiX0}
  (\Pi_\pm \nabla x_0)^\mu
  =  \frac{ 
  \tfrac{1}{2} \left(
  \pm x_0 + \sqrt{(x_0)^2 - \tfrac{4}{[2]^2} (\Xsquare)}
  \right)\delta_0^\mu
  \mp \tfrac{1}{q[2]} x^\mu }
  {\sqrt{(x_0)^2 - \tfrac{4}{[2]^2} (\Xsquare)}} \,.
\end{equation} 

Finally, we can combine the previous result~\eqref{eq:XsquareDer} for the derivative of functions of the four-length $\Xsquare = x_\mu x^\mu$ and \eqref{eq:Lquotient3} for the derivative of functions of the time coordinate $x_0$ into a formula for functions $f = f(\Xsquare, x_0)$ of both variables,
\begin{equation}
\label{eq:Lquotient4}  
\begin{split}  
  \nabla f(\Xsquare,x_0) 
  &= 
  \frac{f(q^2\Xsquare, x_0) - f(\Xsquare, x_0)}
  {q^2\Xsquare - \Xsquare} \,\nabla\Xsquare \\  
  &\quad + 
  \frac{f(q^2\Xsquare, q \tau_+) - f(q^2\Xsquare, x_0)}
  {q\tau_+ - x_0} \Pi_+ \nabla x_0 \\
  &\quad +
  \frac{f(q^2\Xsquare, q \tau_-) - f(q^2\Xsquare, x_0)}
  {q\tau_- - x_0} \Pi_- \nabla x_0 \,.
\end{split}
\end{equation}
This formula can be further simplified to
\begin{equation}
\label{eq:Lquotient5}  
\begin{split}  
  \nabla f(\Xsquare,x_0) 
  = 
  &\frac{f(q^2\Xsquare, q \tau_+) - f(\Xsquare, x_0)}
  {q\tau_+ - x_0} \,\Pi_+ \nabla x_0 \\
  +
  &\frac{f(q^2\Xsquare, q \tau_-) - f(\Xsquare, x_0)}
  {q\tau_- - x_0} \,\Pi_- \nabla x_0 \,.
\end{split}
\end{equation} 
%
%

\subsubsection{Separation of variables}

We have seen that the derivatives of functions of $x_\pm$, $x_\Lind$, and $\Xsquare$ can be written in terms of Jackson $q$-derivatives. Derivatives of functions of the time coordinate, however, do not have this property but depend on the four-length $\Xsquare$, as well. In the language of partial differential equations, the partial derivative $\partial^\mu$ is not separated with respect to the standard time coordinate of quantum Minkowski space. We will now show that there is a remarkable non-linear transformation of coordinates such that the partial derivatives are separated in the new coordinates.

These new coordinates $\xi_\pm = \xi_\pm(\Xsquare,x_0)$ are given by
\begin{equation}
\label{eq:Sep1}  
  \xi_\pm = \tfrac{1}{2} 
  \left( x_0 \pm \sqrt{(x_0)^2 - \tfrac{4}{[2]^2} (\Xsquare)}
  \right) \,.
\end{equation} 
In terms of the new variables $\tau_\pm$ is expressed as
\begin{equation}
\begin{aligned}  
  \tau_+ &= q \xi_+ + q^{-1} \xi_- \\
  \tau_- &= q^{-1} \xi_+ + q \xi_- \,. 
\end{aligned}
\end{equation}
The back transform is given by
\begin{equation}
\begin{aligned}  
  x_0      &= \xi_+ + \xi_- \\
  \Xsquare &= [2]^2 \xi_+ \xi_- \,. 
\end{aligned}
\end{equation}
A function $f=f(\Xsquare, x_0)$ is expressed in terms of the new variables by the transformed function
\begin{equation}
\begin{split}  
  \tilde{f}(\xi_+,\xi_-) 
  &:= f\bigl( \Xsquare(\xi_+,\xi_-), x_0(\xi_+,\xi_-) \bigr) \\
  &= f( [2]^2 \xi_+ \xi_-, \xi_+ + \xi_- ) \,,
\end{split}
\end{equation} 
for which the quotients of Eq.~\eqref{eq:Lquotient5} take the form of Jackson derivatives
\begin{equation}
\label{eq:Lquotient6}  
\begin{aligned}  
  \frac{f(q^2\Xsquare, q \tau_+) - f(\Xsquare, x_0)}
       {q\tau_+ - x_0} 
  &= \frac{\tilde{f}(q^2\xi_+, \xi_-) - \tilde{f}(\xi_+, \xi_-)} 
          {q^2 \xi_+ - \xi_+} \\ 
  \frac{f(q^2\Xsquare, q \tau_+) - f(\Xsquare, x_0)}
       {q\tau_+ - x_0} 
  &= \frac{\tilde{f}(\xi_+, q^2\xi_-) - \tilde{f}(\xi_+, \xi_-)} 
          {q^2 \xi_- - \xi_-} \,.
\end{aligned}
\end{equation}
In order to completely rewrite Eq.~\eqref{eq:Lquotient5} in $q$-derivative form we still have to rewrite the expressions $\Pi_\pm \nabla x_0$. Towards this end, we must understand on what the projection operators $\Pi_\pm$ actually project.

Key to understanding the operators $\Pi_\pm$ is the calculation of the derivatives of the new coordinates viewed as functions of $\Xsquare$ and $x_0$. Inserting the functions $\xi_\pm = \xi_\pm(\Xsquare, x_0)$ into formula~\eqref{eq:Lquotient5} and using the formula 
\begin{equation}
  \sqrt{\tau_\pm^2 - \tfrac{4}{[2]^2} (\Xsquare)}
  = \frac{1}{2} \left( \pm \lambda x_0 
  + [2] \sqrt{(x_0)^2 - \tfrac{4}{[2]^2} (\Xsquare)} \right) \,,
\end{equation} 
we obtain after long calculations the compact result
\begin{equation}
  \nabla \xi_\pm = \Pi_\pm \nabla x_0 \,.
\end{equation}
In this sense $\Pi_\pm$ can be viewed as generalized Jacobian of the coordinate transformation~\eqref{eq:Sep1}. We thus arrive at the main result of this paper:

\begin{Proposition}
\label{th:CentralDeriv}
Let $f = f(\xi_+, \xi_-)$ be a general ordered polynomial function in $\xi_\pm$. Then
\begin{equation}
\label{eq:ADeriv}  
  \nabla f(\xi_+, \xi_-)   
  = \frac{\partial_{q^2} f }
         {\partial_{q^2}\xi_+} \nabla \xi_+
  + \frac{\partial_{q^2} f } 
         {\partial_{q^2}\xi_-} \nabla \xi_- \,.
\end{equation}
\end{Proposition}

\subsubsection{General functions}

Finally, we can combine this result with formulas~\eqref{eq:SpaceDeriv} for derivations of functions of the non-central variables. A general polynomial in $\Mink$ can be written as a sum of products of functions $g =f(\xi_+,\xi_-)$ and $h = h_+(x_+,x_\Lind) + h(x_\Lind, x_-)$. Using the coproduct~\eqref{eq:Pcoprod} we obtain for the derivative of the product
\begin{equation}
\label{eq:General1}  
  \nabla(f g) = (\nabla f) g + (\kappa L \tr f) \nabla g \,,
\end{equation}
where we can derive from Proposition~\ref{th:functionL} the formula
\begin{equation}
\label{eq:General2}
  \kappa L \tr g(\xi_+,\xi_-)
  =  \Pi_+ g(q^2 \xi_+,\xi_-)  + \Pi_- g(\xi_+,q^2\xi_-)  \,.
\end{equation}
Inserting Eq.~\eqref{eq:General2} into \eqref{eq:General1} we have to observe that $\Pi_\pm$ does not commute with the non-central observables. In order to write down the end result, we first define for a general linear combination of the basis $\mathcal{B}$ the matrix valued object
\begin{equation}
\label{eq:PifDef}  
  \delta f := 
  \Pi_+\bigl( f |_{\xi_+ \mapsto q^2\xi_+} - f \bigr) + 
  \Pi_-\bigl( f |_{\xi_- \mapsto q^2\xi_-} - f \bigr) \,.
\end{equation}
This definition does not depend on the ordering because $\xi_+$ and $\xi_-$ are both central. Putting things together we obtain for a general $f \in \mathrm{Span}\,\mathcal{B}$
\begin{equation}
\label{eq:GeneralDeriv}
\begin{split}
  \nabla f 
  &=\frac{\partial_{q^2} f }
         {\partial_{q^2}\xi_+} \nabla \xi_+
  + \frac{\partial_{q^2} f } 
         {\partial_{q^2}\xi_-} \nabla \xi_- \\
  &+	 
    \frac{\partial_{q^2} f }
         {\partial_{q^2}x_+} \nabla x_+
  + \frac{\partial_{q^2} f } 
         {\partial_{q^2} x_\Lind} \nabla x_\Lind
  + \frac{\partial_{q^2} f }
         {\partial_{q^2}x_-} \nabla x_- \\
  &+	 
    \frac{\partial_{q^2} (\delta f) }
         {\partial_{q^2}x_+} \nabla x_+
  + \frac{\partial_{q^2} (\delta f) } 
         {\partial_{q^2} x_\Lind} \nabla x_\Lind
  + \frac{\partial_{q^2} (\delta f) }
         {\partial_{q^2}x_-} \nabla x_- \,.	 
\end{split} 
\end{equation} 
The first two lines look like a generalized chain rule with partial Jackson derivatives as outer derivatives. The partial Jackson derivatives in the last line are to be understood not to act on the projections $\Pi_\pm$ contained in $\delta f$. By definition of $\delta f$ the last line vanishes for $q \rightarrow 1$ and the first two lines reproduce the usual chain rule.

\section{Solutions of quantum wave equations}
\label{sec:Wave}

In order to give an application of the separation of variables and the formulas for derivations, let us now come back to the initial problem of calculating solutions of the quantum Klein-Gordon equation~\eqref{eq:KleinGordon2}. It can be shown that, just as in the commutative case, the space of solutions is generated as representation of the quantum algebra by a momentum eigenstate \cite{Blohmann:2001b}. However, unlike in the commutative case, there is only a finite number of such momentum eigenstates. In the massless case it is the light cone state, defined by
\begin{equation}
\label{eq:MomentumEigen1}  
  p_0 \tr \psi = k \psi \,,\qquad p_3 \tr \psi = k \psi \,,
  \qquad p_\pm \tr \psi = 0 \,,
\end{equation}
where $k \in \mathbb{R}$ is related to the helicity. In the massive case it is the rest state
\begin{equation}
\label{eq:MomentumEigen2}  
  p_0 \tr \psi = m \psi \,,\qquad p_3 \tr \psi = 0  \,,
  \qquad p_\pm \tr \psi = 0 \,,
\end{equation}
where $m \in \mathbb{R}$ is the mass. The difficult part is to calculate these states within the noncommutative differential calculus as solutions of quantum differential equations. All other solutions of the quantum Klein-Gordon equation can then be obtained by quantum Lorentz transformations, which are well known and straightforward to calculate.

\subsection{The massless case}

Let us start with the massless case which turns out to be quite simple. As $q$-differential equation Eq.~\eqref{eq:MomentumEigen1} reads
\begin{equation}
\label{eq:EigenDiff1}  
  \partial^0 \tr \psi = \I k \psi \,,\qquad 
  \partial^3 \tr \psi = -\I k \psi \,, \qquad 
  \partial^\pm \tr \psi = 0 \,,
\end{equation} 
where we have raised the indices with the $q$-Minkowski metric, e.g., $\partial^3 = - \partial_3$. Using the coordinate free notation and the lightcone coordinate $x_\Lind$ to be a function of the light cone coordinate $x_\Lind = x_3 - x_0$, we can write Eq.~\eqref{eq:EigenDiff1} as
\begin{equation}
\label{eq:EigenDiff1b}
  \nabla \psi = - \I k \psi \nabla x_\Lind \,.
\end{equation}
Since $x_\Lind$ is one of the variables separating the noncommutative differential calculus, we can conclude that~\eqref{eq:EigenDiff1b} can be solved by a function $\psi = \psi(x_\Lind)$ of $x_\Lind$ alone. Using the formulas~\eqref{eq:SpaceDeriv} for the derivatives we obtain
\begin{equation}
  \nabla \psi(x_\Lind) = 
  \frac{\partial_{q^2} \psi }{\partial_{q^2}x_\Lind} \nabla x_\Lind
  = - \I k \psi \nabla x_\Lind \,,
\end{equation} 
which is equivalent to
\begin{equation}
\label{eq:qExp1}  
  \frac{\partial_{q^2} \psi }{\partial_{q^2}x_\Lind}
  = - \I k \psi(x_\Lind) \,.
\end{equation}
This equation is solved by the well-known $q$-exponential function
\begin{equation}
\label{eq:qExp2}
  \exp_q(z) \equiv \E_q^z := \sum_{n = 1}^\infty \frac{z^n}{\qNum{n}!} \,,
\end{equation}
where the $q$-factorial is defined in the obvious way as 
\begin{equation}
  \qNum{n}! = \qNum{n}\qNum{n-1}\cdots \qNum{1} \,.
\end{equation}
With the $q$-exponential the general solution of Eq.~\eqref{eq:qExp1} and, hence, of Eq.~\eqref{eq:EigenDiff1} can be written as
\begin{equation}
\label{eq:MasslessSolution}
  \psi = C\, \E_q^{-\I k (x_3 - x_0)} \,,
\end{equation} 
where $C$ is a normalization constant.

\subsection{The massive case} 

Let us now turn to the massive case where calculation of the rest state is more involved. As $q$-differential function the defining eigenvalue equation~\eqref{eq:MomentumEigen2} reads
\begin{equation}
\label{eq:EigenDiff2}  
  \partial^0 \tr \psi = \I m \psi \,, \qquad 
  \partial^A \tr \psi = 0 \,,
\end{equation} 
where the index runs through $A = -,3,+$. In index free notation we can write this as
\begin{equation}
\label{eq:EigenDiff2b}
  \nabla \psi = \I m \psi \nabla x_0 \,.
\end{equation}
Since the time coordinate $x_0$ is not one of the variables which separates the noncommutative differential calculus, we cannot conclude that Eq.~\eqref{eq:EigenDiff2b} can be solved by a function in $x_0$ alone. Since the eigenvalue equation is invariant with respect to quantum rotations we can conclude, though, that $\psi$ must be a function of the scalars with respect to quantum rotations $x_0$ and $\Xsquare$. In order to solve Eq.~\eqref{eq:EigenDiff2b} we have to separate it with the separating coordinates $\xi_-$ and $\xi_+$. Using Eq.~\eqref{eq:ADeriv} and $x_0 = \xi_- + \xi_-$ we obtain
\begin{equation}
\label{eq:EigenDiff2c}
  \nabla \psi(\xi_+,\xi_-) 
  = \frac{\partial_{q^2} \psi }
         {\partial_{q^2}\xi_+} \nabla \xi_+
  + \frac{\partial_{q^2} \psi } 
         {\partial_{q^2}\xi_-} \nabla \xi_- 
  = \I m (\psi \nabla \xi_+ + \psi \nabla \xi_-) \,,
\end{equation}
which can be rewritten as 
\begin{equation}
\label{eq:EigenDiff2d}
  ( \frac{\partial_{q^2} \psi }{\partial_{q^2}\xi_+} 
    - \I m \psi) \nabla \xi_+
 +( \frac{\partial_{q^2} \psi }{\partial_{q^2}\xi_-} 
    - \I m \psi) \nabla \xi_- = 0 \,.
\end{equation}
The differential equation is now separated so we can make the ansatz
\begin{equation}
  \psi(\xi_+,\xi_-) = \psi_+ (\xi_+) \,\psi_-(\xi_-) \,,
\end{equation}
which solves Eq.~\eqref{eq:EigenDiff2d} if
\begin{equation}
\label{eq:mWave4}
    \frac{\partial_{q^2} \psi_+}{\partial_{q^2}\xi_+} = \I m \psi_+ 
   \,,\qquad
    \frac{\partial_{q^2} \psi_-}{\partial_{q^2}\xi_-} = \I m \psi_- \,.
\end{equation}
We see that the separation of variables of the differential calculus by the new coordinates $\xi_\pm$ leads to a complete separation of the wave equation. Eqs.~\eqref{eq:mWave4} are again solved by the $q$-exponential function~\eqref{eq:qExp2}. As end result we obtain
\begin{equation}
\label{eq:MassiveSolution}
\begin{split}  
  \psi &= C \, \E_q^{\I m \xi_+} \E_q^{\I m \xi_-}\\
  &= C \,    
  \E_q^{\frac{\I m}{2} 
  \left( x_0 + \sqrt{(x_0)^2 - \tfrac{4}{[2]^2} (\Xsquare)} \right)}
  \E_q^{\frac{\I m}{2} 
  \left( x_0 - \sqrt{(x_0)^2 - \tfrac{4}{[2]^2} (\Xsquare)} \right)} \,,
\end{split}
\end{equation}
where $C$ is a normalization constant. This rest state is a solution of the massive quantum Klein-Gordon equation~\eqref{eq:KleinGordon2}. All other solutions can be obtained from this one by quantum Lorentz transformations \cite{Bachmaier}.

\section{Conclusion}
\label{sec:Conclusion}

The separation of variables~\eqref{eq:ADeriv} by a nonlinear coordinate transformation~\eqref{eq:Sep1} seems to be the most intriguing result presented in this paper. We do not yet understand on a fundamental level what property of quantum Minkowski space is responsible for the existence of such a transformation. Nor do we know if and how this can be generalized to other quantum spaces.

In definition~\eqref{eq:Sep1} of the separating variables a square root expression appears which, strictly speaking, is not an element of the algebra of Minkowski space proper. In order to make all statements completely rigorous we have to enlarge $\Mink$ by a central square root element $\alpha$ satisfying
\begin{equation}
  \alpha^2 = [2]^2 (x_0)^2 - 4 (\Xsquare) \,. 
\end{equation}
Alternatively, one could think of the coordinates being represented by operators on a Hilbert space as in \cite{Cerchiai:1998}. In this case the square root would be defined through functional calculus for normal operators. In the expansion of solution~\eqref{eq:MassiveSolution} of the quantum Klein-Gordon equation the square root drops out. There, it can be seen as an auxiliary object which makes the notation of a certain generating function more compact. Anyway, it is clear that this mathematical subtlety does not affect the results of this paper.

Finally, we would like to note that expressions for derivatives within the noncommutative calculus can also be derived in a straightforward manner using the recursion relations defined by the commutation relations of momenta and coordinates. While this approach circumvents to some degree the mathematical machinery we have introduced here, it produces lengthy formulas which do not give the structural insight desirable for solving noncommutative differential equations. Even when free field equation can be solved by brute force, using computer algebra for example, the results turn to be quite complicated. But if already something as basic as the wave function of a free particle were described by complicated expressions, some fundamental questions of quantum field theory, such as independence of in and out states, would become very hard to address. In this respect we believe results of the type presented here to be not just a matter of computational convenience but an indispensable requirement for the further development of quantum theory on quantum Minkowski space.

\appendix

\section{Basic definitions and formulas}
\label{sec:AppBasic}

\begin{Definition}
The Hopf $*$-algebra generated by $E$, $F$, $K$, and $K^{-1}$ with relations
\begin{equation}
\begin{gathered}
  KK^{-1}= 1 = K^{-1}K \,,\quad KEK^{-1} = q^2 E \,,\\
  KFK^{-1} = q^{-2}F\,,\quad [E,F]= \lambda^{-1}(K-K^{-1})\,,
\end{gathered}
\end{equation}
Hopf structure
\begin{equation}
\begin{gathered}
  \Delta(E) = E\otimes K + 1\otimes E \,,\quad
  \Delta(F) = F\otimes 1 + K^{-1}\otimes F \,, \\
  \Delta(K) = K \otimes K \,,\quad
  \varepsilon(E)=0=\varepsilon(F) \,,\quad \varepsilon(K)=1\,,\\
  S(E) = -EK^{-1} \,,\quad S(F) = -KF \,,\quad S(K)=K^{-1}\,,
\end{gathered}
\end{equation}
and $*$-structure
\begin{equation}
  E^* = FK \,, F^* = K^{-1}E \,, K^* = K
\end{equation}
is called $\Uqsu$, the $q$-deformation of the enveloping algebra $\Usu$   \cite{Kulish:1983,Sklyanin:1985}.
\end{Definition}
Another useful set of generators is given by angular momentum vector $\{J_A\} = \{ J_-, J_3, J_+\}$ which is a vector operator with respect to the Hopf adjoint action of $\Usu$ on itself \cite{Blohmann:2001b}, 
\begin{equation}
\label{eq:DefJ}
\begin{aligned}
  J_{-} &:= q[2]^{-\frac{1}{2}}KF \,,\\
  J_3   &:= [2]^{-1} (q^{-1}EF-qFE) \,,\\
  J_{+} &:= -[2]^{-\frac{1}{2}}E \,.
\end{aligned}
\end{equation}
The spin-$\frac{1}{2}$ representation of the antipode of this vector operator yields by $\tilde{\sigma}\!_A := - [2] \rho^{\frac{1}{2}}(SJ_A)$ a variant of the $q$-Pauli matrices which we will need here. Explicitly, we obtain
\begin{equation}
\label{eq:sigmaDef}
  \tilde{\sigma}_- = {[2]}^{\frac{1}{2}}
    \begin{pmatrix}
      0 & q^{\frac{1}{2}} \\ 0 & 0
    \end{pmatrix}, 
  \tilde{\sigma}_+ = [2]^{\frac{1}{2}}
    \begin{pmatrix}
      0 & 0 \\ -q^{-\frac{1}{2}} & 0
    \end{pmatrix},
  \tilde{\sigma}_3 =
    \begin{pmatrix} -q^{-1} & 0 \\ 0 & q \end{pmatrix} 
\end{equation}
with respect to the $\{-,+\}$ basis. $\Uqsu$ is quasitriangular with universal $\R$-matrix \cite{Drinfeld:1986}
\begin{equation}
\label{eq:Rmatrix}
  \R = \E^{\h (H\otimes H)/2} \sum_{n=0}^{\infty}
  \E^{\h n(n-1)/2} \frac{(\E^\h- \E^{-\h})^n}{[n]!} (E^n \otimes F^n) \,,
\end{equation}
where $\h := \ln q$. The quantum Lorentz algebra contains a $\SUq^\op$ Hopf sub-algebra \cite{Podles:1990}, the generators of which are given explicitly by
\begin{equation}
\label{eq:boostdef}
\begin{aligned}
  a &:= K^{\frac{1}{2}}\otimes K^{-\frac{1}{2}} \\
  b &:= q^{-\frac{1}{2}}\lambda K^{\frac{1}{2}}
    \otimes K^{-\frac{1}{2}} E \\
  c &:=  -q^{\frac{1}{2}}\lambda F K^{\frac{1}{2}}
    \otimes K^{-\frac{1}{2}} \\
  d &:= K^{-\frac{1}{2}}\otimes K^{\frac{1}{2}}
    - \lambda^2 F K^{\frac{1}{2}}\otimes K^{-\frac{1}{2}} E \,.
\end{aligned}
\end{equation}

\subsubsection*{Acknowledgements}

The work of F.B.\ was supported by a PhD fellowship of the Max-Planck Society. C.B.\ was supported by the European Union with an outgoing international Marie-Curie fellowship under contract MOIF-CT-2005-8559.

\providecommand{\href}[2]{#2}\begingroup\raggedright\endgroup


\begin{thebibliography}{10}

\bibitem{Heisenberg:1938}
W.~Heisenberg, ``\"Uber die in der Theorie der Elementarteilchen auftretende
  universelle L\"ange,'' {\em Ann. Phys.} {\bf 32} (1938) 20--33.

\bibitem{Minwalla:1999}
S.~Minwalla, M.~Van~Raamsdonk, and N.~Seiberg, ``Noncommutative perturbative
  dynamics,'' {\em JHEP} {\bf 02} (2000) 020,
\href{http://www.arXiv.org/abs/hep-th/9912072}{{\tt hep-th/9912072}}.

\bibitem{Matusis:2000}
A.~Matusis, L.~Susskind, and N.~Toumbas, ``The IR/UV connection in the
  non-commutative gauge theories,'' {\em JHEP} {\bf 12} (2000) 002,
\href{http://www.arXiv.org/abs/hep-th/0002075}{{\tt hep-th/0002075}}.

\bibitem{Grosse:2005}
H.~Grosse and R.~Wulkenhaar, ``Renormalisation of $\phi^4$-Theory on
  Non-Commutative $\mathbb{R}^4$ to All Orders,'' {\em Lett. Math. Phys.} {\bf
  71} (2005)
13--26.

\bibitem{Carow-Watamura:1990}
U.~Carow-Watamura, M.~Schlieker, M.~Scholl, and S.~Watamura, ``Tensor
  Representation of the Quantum Group $SL_q(2)$ and Quantum Minkowski Space,''
  {\em Z. Phys.} {\bf C48} (1990)
159--166.

\bibitem{Doplicher:1995}
S.~Doplicher, K.~Fredenhagen, and J.~E. Roberts, ``The Quantum structure of
  space-time at the Planck scale and quantum fields,'' {\em Commun. Math.
  Phys.} {\bf 172} (1995) 187--220,
\href{http://www.arXiv.org/abs/hep-th/0303037}{{\tt hep-th/0303037}}.

\bibitem{Ogievetskii:1992a}
O.~Ogievetskii, W.~B. Schmidke, J.~Wess, and B.~Zumino, ``$q$-Deformed
  Poincar{\'e} algebra,'' {\em Commun. Math. Phys.} {\bf 150} (1992)
495.

\bibitem{Faddeev:1990}
L.~D. Faddeev, N.~Y. Reshetikhin, and L.~A. Takhtajan, ``Quantization of Lie
  Groups and Lie Algebras,'' {\em Leningrad Math. J.} {\bf 1} (1990)
193--225.

\bibitem{Majid:1993}
S.~Majid, ``Braided momentum in the $q$-Poincare group,'' {\em J. Math. Phys.}
  {\bf 34} (1993) 2045--2058,
\href{http://www.arXiv.org/abs/hep-th/9210141}{{\tt hep-th/9210141}}.

\bibitem{Cerchiai:1998}
B.~L. Cerchiai and J.~Wess, ``$q$-Deformed Minkowski Space based on a
  $q$-Lorentz Algebra,'' {\em Eur. Phys. J.} {\bf C5} (1998) 553--566,
\href{http://www.arXiv.org/abs/math.qa/9801104}{{\tt math.qa/9801104}}.

\bibitem{Blohmann:2001a}
C.~Blohmann, ``Spin in the $q$-Deformed Poincar{\'e} Algebra,'' {\em Commun.
  Math. Phys.} {\bf 243} (2003) 329--342,
\href{http://www.arXiv.org/abs/math.qa/0111008}{{\tt math.qa/0111008}}.

\bibitem{Blohmann:2001b}
C.~Blohmann, ``Free $q$-Deformed Relativistic Wave Equations by Representation
  Theory,'' {\em Eur. Phys. J.} {\bf C30} (2003) 435--445,
\href{http://www.arXiv.org/abs/hep-th/0111172}{{\tt hep-th/0111172}}.

\bibitem{Madore:2000b}
J.~Madore, S.~Schraml, P.~Schupp, and J.~Wess, ``Gauge theory on noncommutative
  spaces,'' {\em Eur. Phys. J.} {\bf C16} (2000) 161--167,
\href{http://www.arXiv.org/abs/hep-th/0001203}{{\tt hep-th/0001203}}.

\bibitem{Jurco:2001}
B.~Jurco, L.~Moller, S.~Schraml, P.~Schupp, and J.~Wess, ``Construction of
  non-Abelian gauge theories on noncommutative spaces,'' {\em Eur. Phys. J.}
  {\bf C21} (2001) 383--388,
\href{http://www.arXiv.org/abs/hep-th/0104153}{{\tt hep-th/0104153}}.

\bibitem{Majid}
S.~Majid, {\em Foundations of Quantum Group Theory}.
\newblock Cambridge Univ. Press, 1995.

\bibitem{Fiore:1996}
G.~Fiore and P.~Schupp, ``Identical particles and quantum symmetries,'' {\em
  Nucl. Phys.} {\bf B470} (1996) 211--235,
\href{http://arXiv.org/abs/hep-th/9508047}{{\tt hep-th/9508047}}.

\bibitem{Dobrev:1994}
V.~K. Dobrev, ``New q-Minkowski space-time and q-Maxwell equations hierarchy
  from q-conformal invariance,'' {\em Phys. Lett.} {\bf B341} (1994)
133--138.

\bibitem{Schirrmacher:1992}
A.~Schirrmacher, ``Quantum groups, quantum space-time, and Dirac equation,''.
  Talk given at NATO Advanced Research Workshop on Low Dimensional Topology and
  Quantum Field Theory, Cambridge, England, 6-13 Sep 1992.

\bibitem{Pillin:1994b}
M.~Pillin, ``$q$-deformed relativistic wave equations,'' {\em J. Math. Phys.}
  {\bf 35} (1994) 2804--2817,
\href{http://www.arXiv.org/abs/hep-th/9310097}{{\tt hep-th/9310097}}.

\bibitem{Song:1992}
X.-C. Song, ``Covariant differential calculus on quantum Minkowski space and
  the $q$-analog of Dirac equation,'' {\em Z. Phys.} {\bf C55} (1992)
417--422.

\bibitem{Meyer:1995}
U.~Meyer, ``Wave equations on $q$-Minkowski space,'' {\em Commun. Math. Phys.}
  {\bf 174} (1995) 457--476,
\href{http://www.arXiv.org/abs/hep-th/9404054}{{\tt hep-th/9404054}}.

\bibitem{Podles:1996}
P.~Podles, ``Solutions of Klein-Gordon and Dirac equations on quantum Minkowski
  spaces,'' {\em Commun. Math. Phys.} {\bf 181} (1996) 569--586,
\href{http://www.arXiv.org/abs/q-alg/9510019}{{\tt q-alg/9510019}}.

\bibitem{Woronowicz:1989}
S.~L. Woronowicz, ``Differential calculus on compact matrix pseudogroups
  (quantum groups),'' {\em Commun. Math. Phys.} {\bf 122} (1989)
125--170.

\bibitem{Wess:1991}
J.~Wess and B.~Zumino, ``Covariant Differential Calculus on the Quantum
  Hyperplane,'' {\em Nucl. Phys. Proc. Suppl.} {\bf 18B} (1991)
302--312.

\bibitem{Schmuedgen}
A.~Klymik and K.~Schm\"udgen, {\em Quantum Groups and Their Representations}.
\newblock Springer, 1997.

\bibitem{Bauer:2002}
C.~Bauer and H.~Wachter, ``Operator Representations on Quantum Spaces,'' {\em
  Eur. Phys. J.} {\bf C31} (2003) 261--275,
\href{http://www.arXiv.org/abs/math-ph/0201023}{{\tt math-ph/0201023}}.

\bibitem{Schmidke:1991}
W.~B. Schmidke, J.~Wess, and B.~Zumino, ``A $q$-deformed Lorentz algebra,''
  {\em Z. Phys.} {\bf C52} (1991)
471.

\bibitem{Lorek:1997a}
A.~Lorek, W.~Weich, and J.~Wess, ``Non-commutative Euclidean and Minkowski
  structures,'' {\em Z. Phys.} {\bf C76} (1997)
375.

\bibitem{Rohregger:1999}
M.~Rohregger and J.~Wess, ``$q$-deformed Lorentz-algebra in Minkowski phase
  space,'' {\em Eur. Phys. J.} {\bf C7} (1999), no.~1, 177--183.

\bibitem{Podles:1990}
P.~Podles and S.~L. Woronowicz, ``Quantum deformation of Lorentz group,'' {\em
  Commun. Math. Phys.} {\bf 130} (1990)
381.

\bibitem{Bachmaier}
F.~Bachmaier, {\em The free particle on $q$-Minkowski space}.
\newblock PhD thesis, LMU Munich, Faculty of Physics, 2004.
\newblock
  \href{http://www.arXiv.org/abs/http://nbn-resolving.de/urn:nbn:de:bvb:19-191%
76}{{\tt http://nbn-resolving.de/urn:nbn:de:bvb:19-19176}}.

\bibitem{Kulish:1983}
P.~P. Kulish and N.~Y. Reshetikhin, ``Quantum linear problem for the
  Sine-Gordon equation and higher representations,'' {\em J. Sov. Math.} {\bf
  23} (1983)
2435--2441.

\bibitem{Sklyanin:1985}
E.~K. Sklyanin, ``On an Algebra Generated by Quadratic Relations,'' {\em
  Uspekhi Mat. Nauk} {\bf 40} (1985) 214.

\bibitem{Drinfeld:1986}
V.~G. Drinfeld, ``Quantum Groups,'' in {\em Proceedings of the International
  Congress of Mathematicians}, A.~M. Gleason, ed., pp.~798--820.
\newblock Amer. Math. Soc., 1986.

\end{thebibliography}
\end{document}